\documentclass[11pt]{amsart}
\usepackage{amssymb}
\usepackage{latexsym,bm}
\usepackage{amsmath,amsfonts,amssymb}

\newcommand{\BS}{\mathfrak S}\renewcommand{\l}{\ell}
\newcommand{\Z}{\mathbb Z} \renewcommand{\H}{\mathcal{H}}
 
 \newcommand{\lam}{\lambda}

\newcommand{\ts}{\widetilde{S}}
\newcommand{\eps}{\varepsilon}

 \newcommand{\td}{\widetilde{D}}

\newcommand{\ulam}{{\lam}} \newcommand{\sig}{\sigma}

\DeclareMathOperator{\Ind}{Ind} 
\DeclareMathOperator{\res}{res} 
\DeclareMathOperator{\id}{id} \DeclareMathOperator{\Aut}{Aut}
 \DeclareMathOperator{\rad}{rad}

\DeclareMathOperator{\aff}{aff} \DeclareMathOperator{\HH}{h}
\DeclareMathOperator{\bQ}{Q}
\DeclareMathOperator{\WITH}{with}
\DeclareMathOperator{\AN}{an}
\DeclareMathOperator{\APP}{appendix}
\DeclareMathOperator{\BY}{by}
\DeclareMathOperator{\XY}{Xiaoyi}
\DeclareMathOperator{\CI}{Cui}

\newcommand{\bbQ}{\overrightarrow\bQ} 
 \DeclareMathOperator{\Irr}{Irr}
\newcommand{\ksl}{\widehat{\mathfrak{sl}}}

\newtheorem{prop}{Proposition}
\newtheorem{thm}{Theorem}\newtheorem{cor}{Corollary}
\newtheorem{lem}{Lemma}\newtheorem{dfn}{Definition}
\numberwithin{equation}{section} \numberwithin{prop}{section}
\numberwithin{thm}{section} \numberwithin{lem}{section}
\numberwithin{dfn}{section} \numberwithin{cor}{section}

\begin{document}

\title[\null]{The number of simple modules
for the Hecke algebras of type $G(r,p,n)$\\[.3cm] {\small ($\WITH$\! $\AN$\! $\APP$\! $\BY$\! $\XY$\! $\CI$)}}

\author[\null]{Jun Hu}

%\address{Department of Applied Mathematics,
%Beijing Institute of Technology,
%Beijing 100081 P.R. China}
%\email{junhu303@yahoo.com.cn}

\dedicatory{\normalsize Dedicated to Professor Gus Lehrer on the occasion of his sixtieth birthday\\[5pt]}



\begin{abstract}
We derive a parameterization of simple modules for the cyclotomic
Hecke algebras of type $G(r,p,n)$ over field of any characteristic coprime to $p$.
We give explicit formulas for the number of simple
modules over these cyclotomic Hecke algebras.
\end{abstract}

\maketitle

\section{Introduction}

Let $r$, $p$, $d$ and $n$ be positive integers such that $pd=r$.
Let $K$ be a field such that $K$ contains a primitive $p$-th root
of unity $\eps$. Let $x_1,\cdots,x_d$ be invertible elements in
$K$. Let $q\neq 1$ be an
invertible element in $K$. Let $\H_K(r,n)$ be the unital
$K$-algebra with generators $T_0,T_1,\cdots,T_{n-1}$ and relations
$$\begin{aligned}
&(T_{0}^{p}-x_1^p)(T_{0}^{p}-x_2^p)\cdots (T_{0}^{p}-x_d^p)=0,\\
&T_0T_1T_0T_1=T_1T_0T_1T_0,\\
&(T_i+1)(T_i-q)=0,\quad\text{for $1\leq i\leq n-1$,}\\
&T_iT_{i+1}T_i=T_{i+1}T_{i}T_{i+1},\quad\text{for $1\leq i\leq n-2$,}\\
&T_iT_j=T_jT_i,\quad\text{for $0\leq i<j-1\leq n-2$.}\end{aligned}
$$
Let $\H_{K}(r,p,n)$ be the subalgebra of $\H_{K}(r,n)$ generated
by the elements $T_0^{p},\, T_{u}:=T_{0}^{-1}T_{1}T_{0},\,
T_{1},\,T_2,\,\cdots,\,T_{n-1}$. This algebra is called the
cyclotomic Hecke algebra of type $G(r,p,n)$, which was introduced in \cite{A3}, \cite{AK} and
\cite{BM}. It includes Hecke algebras of type $A$, type $B$ and type $D$ as special cases. The algebra $\H_{K}(r,1,n)$ is called the Ariki--Koike algebra. These algebras are conjecturely related to Lusztig's induced characters in the modular representation of finite reductive groups over field of non-defining characteristic (see \cite{BM}).
\smallskip

The representation of Ariki--Koike algebra (e.g., $\H_K(r,n)$) is
well understood by the work of \cite{A1}, \cite{A2}, \cite{DJM} and
\cite{DM}. Let $\mathcal{P}_n$ be the set of
$r$-multipartitions of $n$. Let $\bbQ:=(Q_1,\cdots,Q_r)$ be a fixed
arbitrary permutation of $$ \Bigl(\underbrace{x_1,x_1\eps,\cdots,
x_1\eps^{p-1}}_{\text{$p$ terms}},\cdots,\cdots,
\underbrace{x_d,x_d\eps,\cdots,x_d\eps^{p-1}}_{\text{$p$ terms}}
\Bigr).
$$
We use $\bQ$ to denote the underlying unordered multiset (allowing
repetitions) of $\bbQ$. For any $\lam\in\mathcal{P}_n$, let
$\ts^{\ulam}_{\bbQ}$ be the Specht module defined in \cite{DJM}.
There is a naturally defined bilinear form $\langle,\rangle$ on
$\ts^{\ulam}_{\bbQ}$. Let
$\td^{\lam}_{\bbQ}=\ts^{\ulam}_{\bbQ}/\rad\langle,\rangle$. By
\cite{DJM}, the set
$\bigl\{\td^{\ulam}_{\bbQ}\bigm|\text{$\ulam\in\mathcal{P}_n$,
$\td_{\bbQ}^{\ulam}\neq 0$}\bigr\}$ forms a complete set of pairwise
non-isomorphic simple $\H_{K}(r,n)$-modules. By \cite{A2} and
\cite{DM}, $\td_{\bbQ}^{\lam}\neq 0$ if and only if $\lam$ is a
Kleshchev $r$-multipartition of $n$ with respect to $(q,\bbQ)$.
\medskip

When $q\neq 1$ is a root of unity, Jacon gives in \cite{Ja} another
parameterization of simple $\H_K(r,n)$-modules via FLOTW
$r$-partitions. As an application, a parameterization of simple
$\H_K(r,p,n)$-modules is obtained in \cite{GJ}. The parameterization
results in both \cite{GJ} and \cite{Ja} are valid only when
$K=\mathbb{C}$ (the complex number field). In \cite{Hu3} and
\cite{Hu5}, using a different approach, we obtain a parameterization
of simple $\H_K(p,p,n)$-modules which is valid over field of any
characteristic coprime to $p$, and we give explicit formula for the
number of simple modules of $\H_K(p,p,n)$. In this paper, combining
the results in \cite{GJ} with the results and ideas in \cite{Hu5},
we derive a parameterization as well as explicit formula for the
number of simple $\H_K(r,p,n)$-modules which is valid over field of
any characteristic coprime to $p$. These results generalize the
earlier results in \cite{Ge}, \cite{GJ}, \cite{Hu1}, \cite{Hu2},
\cite{Hu3}, \cite{Hu4}, \cite{Hu5} and \cite{P}, and was already announced in \cite{Hu6}. 
At the end of this
paper there is an appendix given by Xiaoyi Cui who fixes a gap in
the proof of \cite[(2.2)]{Gen}. We remark that the latter result is
crucial to both the present paper and the paper \cite{GJ}.
\medskip

Throughout this paper, $q\neq 1$ is an invertible
element in $K$. Let $e$ be the
smallest positive integer such that $1+q+q^2+\cdots+q^{e-1}=0$ in
$K$; or $\infty$ if no such positive integer exists. We fix elements
$z_1,\cdots,z_s\in K^{\times}$, such that $z_iz_{j}^{-1}\notin
q^{\Z}$, $\forall\,i\neq j$, and for each $1\leq i\leq r$, $Q_i\in
z_jq^{\Z}$ for some $1\leq j\leq s$.

\bigskip\bigskip
\section{Kleshchev $r$-multipartitions and Kleshchev's good lattice}

Let $\lam$ be an $r$-multipartition. The diagram of $\lam$ is the set
$$[\lam]=\bigl\{(i,j,s)\bigm|1\le j\le\lambda^{(s)}_i\text{ for }1\le s\le r\bigr\}.$$
The elements of $[\lam]$ are called the nodes of $\lam$. Given any two nodes
$\gamma=(a,b,c), \gamma'=(a',b',c')$ of $\lam$, say that $\gamma$
is {\it below} $\gamma'$, or $\gamma'$ is {\it above} $\gamma$
with respect to the {\it Kleshchev order}, if either $c>c'$ or
$c=c'$ and $a>a'$. With respect to the
$(r+1)$-tuple $(q,Q_1,\cdots,Q_r)$, the residue of a node $\gamma=(a,b,c)$
is defined to be $\res(\gamma):=Q_cq^{b-a}\in K$. We call $\gamma$ a $\res(\gamma)$-node.
The node $\gamma=(a,\lam^{(c)}_a,c)$ is called a removable node of $\lam$ if $\lam^{(c)}_a>\lam^{(c)}_{a+1}$.
In that case, $\lam\setminus\{\gamma\}$ is again an $r$-multipartition, and we call $\gamma$ an addable node of
$\lam\setminus\{\gamma\}$. For a fixed residue $x\in K$, say that a
removable $x$-node $\gamma$ of $\lam$ is a normal $x$-node, if
whenever $\eta$ is an addable $x$-node of $\lam$ which is below
$\gamma$, there are more removable $x$-node of $\lam$ between
$\eta$ and $\gamma$ than there are addable $x$-nodes. If $\gamma$
is the highest normal $x$-node of $\lam$, we say that $\gamma$ is
a good $x$-node. If $\lam$ is obtained from $\mu$ by removing a
good $x$-node of $\mu$, we write that
$\lam\overset{x}{\twoheadrightarrow}\mu$.

\begin{dfn} {\rm (\cite{AM}, \cite{MM})} \label{df21} Suppose $n\geq 0$. The set
$\mathcal{K}_n$ of Kleshchev $r$ multipartitions of $n$ with respect to
$(q,Q_1,\cdots,Q_r)$ is defined inductively as follows:
\smallskip

(1)
$\mathcal{K}_0:=\Bigl\{\underline{\emptyset}:=\bigl(\underbrace{\emptyset,\cdots,
\emptyset}_{\text{$r$ copies}}\bigl)\Bigr\}$;\smallskip

(2)
$\mathcal{K}_{n+1}:=\Bigl\{\mu\in\mathcal{P}_{n+1}\Bigm|\text{$\lam\overset{x}
{\twoheadrightarrow}\mu$ for some $\lam\in\mathcal{K}_n$ and some
$x\in K$}\Bigr\}$.
\end{dfn}

Let $\mathcal{K}:=\sqcup_{n\geq 0}\mathcal{K}_n$. The {\it Kleshchev's good lattice} with respect to
$(q,Q_1,\cdots,$ $Q_r)$ is, by definition, the infinite graph whose
vertices are the Kleshchev $r$ multipartitions with respect to
$(q,Q_1,\cdots,Q_r)$ and whose arrows are given by $\text{$\lam\overset{x}{\twoheadrightarrow}\mu$}$.
For any $\lam\in\mathcal{K}$ and any $1\leq i\leq s, 0\leq j\leq e-1$, we define $$\begin{aligned}
\widetilde{f}_{i,j}\lam:&=\begin{cases} \lam\cup\{\gamma\}, &\text{if $\gamma$ is a good $(z_iq^{j})$-node of $\lam\cup\{\gamma\}$;}\\
0, &\text{otherwise},
\end{cases}\\
\widetilde{e}_{i,j}\lam:&=\begin{cases} \lam\setminus\{\gamma\}, &\text{if $\gamma$ is a good $(z_iq^{j})$-node of $\lam$;}\\
0, &\text{otherwise}.
\end{cases}
\end{aligned}
$$
By \cite{A2} and \cite{DM}, for any $\lam\in\mathcal{P}_n$,
$\td_{\bbQ}^{\lam}\neq 0$ if and only if $\lam\in\mathcal{K}_n$.

\begin{dfn} Let $\tau$ be the $K$-algebra automorphism of $\H_{K}(r,n)$ which is defined on generators by
$\tau(T_1)=T_0^{-1}T_1T_0, \tau(T_i)=T_i$, for any $i\neq 1$. Let $\sig$ be the nontrivial $K$-algebra automorphism of $\H_{K}(r,n)$
which is defined on generators by $\sig(T_0)=\eps T_0, \sig(T_i)=T_i$, for any $1\leq i\leq n-1$.
\end{dfn}

Note that if $M$ is a simple $\H_{K}(r,n)$-module, then $M^{\sigma}$ is again a simple $\H_{K}(r,n)$-module.

\begin{dfn} Let $\HH$ be the
automorphism of $\mathcal{K}_n$ which is defined by $\bigl(\td_{\bbQ}^{\lam}\bigr)^{\sig}\cong\td_{\bbQ}^{\HH(\lam)}$.
\end{dfn}

Clearly, $\HH^p=\id$. In particular, we get an action of the
cyclic group $C_p$ on $\mathcal{K}_n$ given as follows:
$$\td_{\bbQ}^{\sig^{k}\cdot\lam}\cong
(\td_{\bbQ}^{\lam})^{\sig^k},\quad\,\forall\,k\in\mathbb{Z}.$$ Let
$\sim_{\sig}$ be the corresponding equivalence relation on
${\mathcal{K}}_n$. That is, $\lam\sim_{\sig}\mu$ if and only if
$\lam=g\cdot\mu$ for some $g\in C_p$. For each
$\lam\in{\mathcal{K}}_n/{\sim_{\sig}}$, let ${C}_{\lam}$ be the
stabilizer of $\lam$ in $C_p$. The following results are basically
followed from \cite[(2.2)]{Gen} and \cite[Lemma 2.2]{GJ} (see
\cite[(5.4),(5.5),(5.6)]{Hu3} for an independent proof in the case
where $r=p$). Unfortunately, the proof of \cite[(2.2)]{Gen} given
there contains a gap (as noted in \cite{Ja1}). That is, in the 10th line of Page 527, Genet's
claim about the determinant of the representing matrix is generally
false. Since the result \cite[(2.2)]{Gen} is crucial for both the
present paper and the paper \cite{GJ}, we include at the end of this
paper an appendix given by Xiaoyi Cui who fixes the gap.

\addtocounter{lem}{3}
\begin{lem} \label{lm24} Suppose that $\H_{K}(r,p,n)$ is split over $K$.
\smallskip

1) Let $\td_{\bbQ}^{\lam}$ be any given irreducible
$\H_{K}(r,n)$-module and $D$ be an irreducible
$\H_{K}(r,p,n)$-submodule of $\td_{\bbQ}^{\lam}$. Let $d_0$ be the
smallest positive integer such that $D\cong \bigl(D\bigr)^{\tau^{d_0}}$. Then
$1\leq d_0\leq p$, and $k:=p/d_0$ is the smallest positive integer such
that $\td_{\bbQ}^{\lam}\cong(\td_{\bbQ}^{\lam})^{\sig^k}$, and $$
\td_{\bbQ}^{\lam}\downarrow_{\H_{K}(r,p,n)}\cong D\oplus
D^{\tau}\oplus\cdots\oplus \bigl(D\bigr)^{\tau^{d_0-1}}. $$

2) The set $\Bigl\{D^{\lam,0},D^{\lam,1},\cdots,D^{\lam,|{C}_{\lam}|-1}\Bigm|
\lam\in{\mathcal{K}}_n/{\sim_{\sig}}\Bigr\} $ forms a complete set of pairwise
non-isomorphic simple $\H_{K}(r,p,n)$-modules, where for each
$\lam\in{\mathcal{K}}_n/{\sim_{\sig}}$, $D^{\lam,0}$ is an irreducible $\H_K(r,p,n)$
submodule of $\td_{\bbQ}^{\lam}$, and
$D^{\lam,i}=\bigl(D^{\lam,0}\bigr)^{\tau^i}$ for
$i=0,1,\cdots,|{C}_{\lam}|-1$.
\end{lem}

\begin{proof} By \cite{Gen} and \cite{GJ}, the results in this lemma hold if $K$ is the complex number field.
Furthermore, it is easy to see that all the arguments in
\cite[(2.2)]{Gen} and \cite[Lemma 2.2]{GJ} are actually valid for any
algebraically closed field $K$ of characteristic coprime to $p$. As
a direct consequence of \cite[Lemma 2.2]{GJ} and Frobenius reciprocity,
the statements in this lemma are valid whenever $K$ is an
algebraically closed field of characteristic coprime to $p$. Now
using the fact (see \cite{GL}) that every simple module for the
algebra $\H_{K}(r,n)$ is always absolutely simple, it follows that
these statements remain valid whenever $\H_{K}(r,p,n)$ is split over
$K$.
\end{proof}

Therefore, the problem on classifying simple
$\H_{K}(r,p,n)$-modules reduces to the problem of determining the
automorphism $\HH$.

\bigskip\bigskip
\section{FLOTW $r$-partitions and FLOTW's good lattice}

Throughout this section, we assume that $e<\infty$.\smallskip

For each integer $c$ with $1\leq c\leq r$, we fix an integer $0\leq v_c\leq e-1$ such that
$Q_c=z_iq^{v_c}$ for some integers $i$ with $1\leq i\leq s$. Throughout this section, we make the following assumption: $$\begin{matrix}
\text{the order on $\bbQ$ are
chosen such that, whenever $Q_{i_1},Q_{i_2},\cdots,Q_{i_s}$ is a $q$}\\
\text{orbit in $\bQ$, where $i_1<i_2<\cdots<i_s$, we have
$0\leq v_{i_1}\leq v_{i_2}\leq\cdots\leq v_{i_s}<e$.} \end{matrix}
$$

Let $\lam=(\lam^{(1)},\cdots,\lam^{(r)})$ be an $r$-multipartition,
or equivalently (see \cite{Ja}), an $r$-partition. Given any two
nodes $\gamma=(a,b,c), \gamma'=(a',b',c')$ of $\lam$, say that
$\gamma$ is {\it below} $\gamma'$, or $\gamma'$ is {\it above}
$\gamma$ with respect to the {\it FLOTW order}, if whenever $Q_c,
Q_{c'}$ are in a single $q$-orbit, then either
$b-a+v_c>b'-a'+v_{c'}$ or $b-a+v_c=b'-a'+v_{c'}$ and $c<c'$. Note
that the FLOTW order does depend on the choice of the elements
$\{z_1,\cdots,z_s\}$ (which we have fixed at the end of Section 1).
In a similar way as before (see Section 2), we have the notions of
normal $x$-nodes and good $x$-nodes with respect to the FLOTW order.
If $\lam$ is obtained from $\mu$ by removing a good $x$-node of
$\mu$, we write that $\lam\underset{x}{\twoheadrightarrow}\mu$. If
$Q_{1_1}, Q_{i_2}, \cdots, Q_{i_s}$ form a $q$-obit in $\bbQ$, then
we call the multipartition
$(\lam^{(i_1)},\lam^{(i_2)},\cdots,\lam^{(i_s)})$ the restriction of
$\lam$ to that $q$-orbit.

\begin{dfn} {\rm (\cite{FLO})} Suppose $n\geq 0$. Let
$\lam=(\lam^{(1)},\lam^{(2)},\cdots,\lam^{(r)})\in\mathcal{P}_n$.
If $\bQ$ is a single $q$-orbit, then $\lam$ is a FLOTW
$r$-partition of $n$ with respect to $(q,\bbQ)$ if and only if:
\begin{enumerate}
\item[(1)] for all $1\leq j\leq r-1$ and $i=1,2,\cdots$, we have $$
\lam_i^{(j)}\geq\lam_{i+v_{j+1}-v_{j}}^{(j+1)},\quad
\lam_i^{(r)}\geq\lam_{i+e+v_{1}-v_{r}}^{(1)};
$$
\item[(2)] for any $k\geq 0$, among the
residues appearing at the right ends of the length $k$ rows of
$\ulam$, at least one element of $\{z_1,z_1q,\cdots,z_1q^{e-1}\}$
does not occur.
\end{enumerate}
In general, if $\bbQ$ is a disjoint union of several $q$-orbits,
then $\lam$ is a FLOTW $r$-partition of $n$ with respect to
$(q,\bbQ)$ if and only if with respect to each $q$-orbit of $\bbQ$,
the restriction of $\lam$ to that $q$-orbit satisfies the above two
conditions.
\end{dfn}

By \cite{FLO}, one can also give a recursive definition (like
(\ref{df21})) of FLOTW $r$-partition by using the procedure of
adding good nodes. Note that, at the moment, we do not have a
non-recursive definition for Kleshchev $r$-multipartition except for
$r\leq 2$ (see \cite{AJ}, \cite{DJ}).
\smallskip

Let $\mathcal{F}_n$ be the set of all the FLOTW $r$-partitions of
$n$ with respect to $(q,Q_1,\cdots,Q_r)$. Let
$\mathcal{F}:=\sqcup_{n\geq 0}\mathcal{F}_n$. The {\it FLOTW's good
lattice} (w.r.t. $(q,Q_1,\cdots,Q_r)$) is, by definition, the
infinite graph whose vertices are the FLOTW $r$-partitions with
respect to $(q,Q_1,\cdots,Q_r)$ and whose arrows are given by
$\text{$\lam\underset{x}{\twoheadrightarrow}\mu$}$. For any
$\lam\in\mathcal{F}$ and any $1\leq i\leq s, 0\leq j\leq e-1$, we
define $$\begin{aligned}
\widetilde{f}_{i,j}\circ\lam:&=\begin{cases} \lam\cup\{\gamma\}, &\text{if $\gamma$ is a good $(z_iq^{j})$-node of $\lam\cup\{\gamma\}$;}\\
0, &\text{otherwise},
\end{cases}\\
\widetilde{e}_{i,j}\circ\lam:&=\begin{cases} \lam\setminus\{\gamma\}, &\text{if $\gamma$ is a good $(z_iq^{j})$-node of $\lam$;}\\
0, &\text{otherwise}.
\end{cases}
\end{aligned}
$$

When $s=1$, i.e., the parameters $\{Q_1,\cdots,Q_r\}$ are in a single $q$-orbit, both the Kleshchev's good lattice and FLOTW's good lattice provide realizations of the crystal graph of certain irreducible integrable highest weight module over the quantum affine algebra $U'_v(\ksl_e)$.

\addtocounter{lem}{1}
\begin{lem}{\rm (\cite[(4.1)]{Ja})} \label{lm32} There is a unique bijection $\kappa: \mathcal{K}\sqcup\{0\}\rightarrow\mathcal{F}\sqcup\{0\}$, such that,
$\kappa(0)=0$, $\kappa(\underline{\emptyset})=\underline{\emptyset}$ and for any $\lam\in\mathcal{K}$, and any $1\leq i\leq s, 0\leq j\leq e-1$, $$
\kappa\Bigl(\widetilde{f}_{i,j}\lam\Bigr)=\widetilde{f}_{i,j}\circ\kappa\bigl(\lam\bigr),\quad
\kappa\Bigl(\widetilde{e}_{i,j}\lam\Bigr)=\widetilde{e}_{i,j}\circ\kappa\bigl(\lam\bigr),
$$
and for any $n\geq 0$ and any $\lam\in\mathcal{K}_n$, the following identity holds in the Grothendieck group of finite dimensional $\H_{K}(r,n)$-modules, $$
\Bigl[\ts_{\bbQ}^{\kappa(\lam)}\Bigr]=\Bigl[\td_{\bbQ}^{\lam}\Bigr]+\sum_{\substack{\mu\in\mathcal{K}_n\\
a(\kappa(\mu))<a(\kappa(\lam))}}d_{\kappa(\lam),\mu}
\Bigl[\td_{\bbQ}^{\mu}\Bigr],$$ where $\lam\in\mathcal{K}_n$,
$d_{\kappa(\lam),\mu}\in\mathbb{Z}^{\geq 0}$ and $a(?)$ is the
$a$-function defined in \cite[(2.4.6)]{GJ}.
\end{lem}

\bigskip\bigskip
\section{A description of $\HH$}

The following are the first two main results in this paper, which
yield (by Lemma \ref{lm24}) a parameterization of simple modules
for the cyclotomic Hecke algebras of type $G(r,p,n)$ over field of
any characteristic coprime to $p$.

\begin{thm} \label{main1} The automorphism $\HH$ does not depend
on the choice of the base field $K$ as long as $\H_K(r,p,n)$ is
split over $K$ and $K$ contains a primitive $p$-th root of unity.
\end{thm}

\begin{proof}  This is proved by using the same argument as in the
Appendix of \cite{Hu4}.
\end{proof}

\begin{thm}\label{main2} Let $\lam\in\mathcal{K}_n$ be a Kleshchev $r$-multipartition of
$n$ with respect to $(q,Q_1,\cdots,Q_r)$. Then, if
$\underline{\emptyset}\overset{r_1}{\twoheadrightarrow}\cdot
\overset{r_2}{\twoheadrightarrow}\cdot \cdots\cdots
\overset{r_n}{\twoheadrightarrow}\lam $ is a path from
$\underline{\emptyset}$ to $\lam$ in Kleshchev's good lattice with
respect to $(q,Q_1,\cdots,$ $Q_r)$, then the sequence
$$\underline{\emptyset}\overset{\eps r_1 }{\twoheadrightarrow}\cdot
\overset{\eps r_2 }{\twoheadrightarrow}\cdot \cdots\cdots
\overset{\eps r_n }{\twoheadrightarrow}\HH(\lam) $$ also defines a
path in Kleshchev's good lattice with respect to
$(q,Q_1,\cdots,Q_r)$, and it connects $\underline{\emptyset}$ to
$\HH(\lam)$.\footnote{One can compare this theorem with
\cite[(1.5)]{Hu2}, \cite[(3.11)]{Hu5} and \cite[Theorem
7.1]{LLT}.}
\end{thm}

Note that the definition of $\HH$ depends on the chosen order on $\bQ$. To stress this point, we had better use the
notation $\HH_{\bbQ}$ instead of $\HH$. Before giving the proof of Theorem \ref{main2}, we note the following result.

\addtocounter{lem}{2}
\begin{lem} Let $\bbQ''=(Q''_1,\cdots,Q''_r)$ be an arbitrary permutation of $\bbQ$. Then Theorem
\ref{main2} is valid for $\HH_{\bbQ}$ if and only if it is valid for $\HH_{\bbQ''}$
\end{lem}

\begin{proof} Let $\mathcal{K}_{\bbQ}$ (resp., $\mathcal{K}_{\bbQ''}$) be the set of Kleshchev $r$-multipartitions with respect to $(q,\bbQ)$ (resp., with respect to $(q,\bbQ'')$).
There is a bijection $\theta$ from $\mathcal{K}_{\bbQ}$ onto
$\mathcal{K}_{\bbQ''}$, such that $\td_{\bbQ}^{\lam}\cong\td_{\bbQ''}^{\theta(\lam)}$,
for any $n\geq 0$ and any $\lam\in\mathcal{K}_{\bbQ}\cap\mathcal{P}_n$. Note that $$
\td_{\bbQ}^{\HH_{\bbQ}(\lam)}=\Bigl(\td_{\bbQ}^{\lam}\Bigr)^{\sigma}\cong\Bigl(\td_{\bbQ''}^{\theta(\lam)}\Bigr)^{\sigma}= \td_{\bbQ''}^{\HH_{\bbQ''}(\theta(\lam))}.$$
It follows that $\theta\bigl(\HH_{\bbQ}(\lam)\bigr)=\HH_{\bbQ''}\bigl(\theta(\lam)\bigr)$.

Since the operator $\widetilde{f}_{i,j}$ can also be defined as
taking socle of the $j$-restriction of module (see \cite{A4}), it
follows that the bijection $\theta$ satisfies $$
\theta\Bigl(\widetilde{f}_{i,j}\lam\Bigr)=\widetilde{f}_{i,j}\theta\bigl(\lam\bigr),
$$
for any $\lam\in \mathcal{K}_{\bbQ}$ and any $1\leq i\leq s, 0\leq j\leq e-1$, from which the lemma follows at once.
\end{proof}

The above lemma allows us to feel free to choose appropriate order
on $\bQ$. This remaining part of this section is devoted to the
proof of Theorem \ref{main2}.\smallskip

First, we need to classify the $q$-orbits in $\bbQ$ (cf.
\cite[Section 2.4.2]{GJ}). Two elements $Q_i, Q_j$ are said to be in
the same $q$-orbit if $Q_i\in Q_jq^{\mathbb{Z}}$; while two elements
$Q_i, Q_j$ are said to be in the same $(\eps,q)$-orbit if $Q_i\in
Q_j\eps^{\mathbb{Z}}q^{\mathbb{Z}}$. Note that for any $1\leq
i,j\leq d$, $x_i\in x_j\eps^{\mathbb{Z}}q^{\mathbb{Z}}$ if and only
if $x_i\eps^a\in (x_j\eps^b)\eps^{\mathbb{Z}}q^{\mathbb{Z}}$ for
some (and hence any) $0\leq a,b\leq p-1$. Therefore, we can split
$\bQ$ into a disjoint union of $d'$ subsets (for some integer $1\leq
d'\leq d$):
$$ \bQ=\bQ^{[1]}\bigsqcup \bQ^{[2]}\bigsqcup\cdots\bigsqcup
\bQ^{[d']},
$$
such that two elements $Q_i, Q_j$ are in the same $(\eps,q)$-orbit
if and only if they belong to the same subset $\bQ^{[i]}$ for some
$1\leq i\leq d'$. Then for each $1\leq i\leq d'$,
$|\bQ^{[i]}|=pd_i$ for some integer $d_i$. Without loss of
generality, we can assume that $$\begin{aligned}
\bQ^{[1]}&=\Bigl\{x_{j}\eps^{b}\Bigm|1\leq j\leq d_{1}, 0\leq b\leq p-1\Bigr\},\\
\bQ^{[2]}&=\Bigl\{x_{j}\eps^{b}\Bigm|d_1+1\leq j\leq d_1+d_{2}, 0\leq b\leq p-1\Bigr\},\\
&\vdots\\
\bQ^{[d']}&=\Bigl\{x_{j}\eps^{b}\Bigm|d-d_{d'}+1\leq j\leq d, 0\leq
b\leq p-1\Bigr\}.
\end{aligned}
$$
Note that for each $\bQ^{[i]}$ and each positive integer $n_i$,
one can naturally associates with a cyclotomic Hecke algebra of
type $G(pd_i,p,n_i)$ with parameters $(q,\bQ^{[i]})$. Let
$\bbQ^{[i]}$ be an ordered $pd_i$-tuple obtained by fixing an
order on $\bQ^{[i]}$, then we have an automorphism
$\HH_{\bbQ^{[i]}}$ on the set of Kleshchev $pd_i$-multipartitions
of $n_i$ with respect to $(q,\bbQ^{[i]})$ which is defined in the
same way as in the case of $G(r,p,n)$. By some abuse of notation,
we use $\bbQ$ to denote
$\bigl(\bbQ^{[1]},\cdots,\bbQ^{[d']}\bigr)$, the concatenation of
ordered tuples. Let $\lam\in\mathcal{K}_n$ be a Kleshchev
$r$-multipartition of $n$ with respect to $(q,\bbQ)$. For each
$1\leq i\leq d'$, let $$\begin{aligned}
\lam^{[i]}:&=\Bigl(\lam^{(pd_1+\cdots+pd_{i-1}+1)},\cdots,\lam^{(pd_1+\cdots+pd_{i})}\Bigr),\\
n_i:&=|\lam^{[i]}|. \end{aligned}$$

\begin{lem} \label{lm44} With the notations as above, we have that $$
\HH_{\bbQ}(\lam)=\Bigl(\HH_{\bbQ^{[1]}}(\lam^{[1]}),\cdots,
\HH_{\bbQ^{[d']}}(\lam^{[d']})\Bigr).
$$
In particular, the Kleshchev $r$-multipartition (with respect to $(q,\bbQ)$) $\HH_{\bbQ}(\lam)$ is as described in Theorem \ref{main2} if and only if for each $1\leq i\leq d'$, the Kleshchev $pd_i$-multipartition (with respect to $(q,\bbQ^{[i]})$) $\HH_{\bbQ^{[i]}}(\lam^{[i]})$ is as described in Theorem \ref{main2}.
\end{lem}

\begin{proof} For each integer $n\geq 0$, let $\H_n^{\aff}$ be the affine Hecke algebra of size $n$ as defined in \cite[Definition 5.1]{Hu5}. By \cite[Corollary 5.6]{Hu5}, we know that $$
D^{\lam}_{\bbQ}\cong\Ind_{\H_{n_1}^{\aff}\otimes\cdots\otimes\H_{n_{d'}}^{\aff}}^{\H_n^{\aff}}\Bigl(
D^{\lam^{[1]}}_{\bbQ_1}\otimes\cdots\otimes
D^{\lam^{[d']}}_{\bbQ_{d'}}\Bigr).
$$
By \cite[Lemma 2.4]{Hu5}, we know that $$
\Bigl(\td_{\bbQ}^{\lam}\Bigr)^{\sigma}\cong\td_{(\eps\!\bbQ^{[1]},\eps\!\bbQ^{[2]},\cdots,\eps\!\bbQ^{[d']})}^{\lam}.
$$
Therefore, $$\begin{aligned}
&\quad\,\,\Bigl(\td_{\bbQ}^{\lam}\Bigr)^{\sigma}\cong\td_{(\eps\!\bbQ^{[1]},\eps\!\bbQ^{[2]},\cdots,\eps\!\bbQ^{[d']})}^{\lam}\\
&\cong\Ind_{\H_{n_1}^{\aff}\otimes\H_{n_{2}}^{\aff}\otimes\cdots\otimes\H_{n_{d'}}^{\aff}}^{\H_n^{\aff}}\Bigl(
D^{\lam^{[1]}}_{\eps\!\bbQ^{[1]}}\otimes
D^{\lam^{[2]}}_{\eps\!\bbQ^{[2]}}\otimes\cdots\otimes
D^{\lam^{[d']}}_{\eps\!\bbQ^{[d']}}\Bigr)\\
&\cong
\Ind_{\H_{n_1}^{\aff}\otimes\H_{n_{2}}^{\aff}\otimes\cdots\otimes\H_{n_{d'}}^{\aff}}^{\H_n^{\aff}}\Bigl(
D^{\HH_{\bbQ^{[1]}}(\lam^{[1]})}_{\bbQ^{[1]}}\otimes
D^{\HH_{\bbQ^{[2]}}(\lam^{[2]})}_{\bbQ^{[2]}}\otimes\cdots\otimes
D^{\HH_{\bbQ^{[d']}}(\lam^{[d']})}_{\bbQ^{[d']}}\Bigr)\\
&\cong \td_{\bbQ}^{\bigl(\HH_{\bbQ^{[1]}}(\lam^{[1]}),\HH_{\bbQ^{[2]}}(\lam^{[2]}),\cdots,
\HH_{\bbQ^{[d']}}(\lam^{[d']})\bigr)},\end{aligned}
$$
from which the lemma follows at once.
\end{proof}

The above lemma allows us to assume without loss of generality that
all the elements $Q_1,\cdots,Q_r$ are in a single $(\eps,q)$-orbit.
{\it Henceforth, we assume that all the parameters in $\bbQ$ are in
a single $(\eps,q)$-orbit.} Recall that we have assumed that $q\neq
1$ from the very beginning. The following two lemmas are useful in
our discussion.

\begin{lem} Let $0\neq a\in
K$. Let $a\!\bQ=\{aQ_1,\cdots,aQ_r\}$. Let $\sigma_a$ be the
isomorphism from $\H_{r,n}^{K}(q,a\!\bQ)$ onto $\H_{r,n}^{K}(q,\bQ)$
which is defined on generators by $\sigma_a(T_0)=aT_0$ and
$\sigma_a(T_i)=T_i$ for any $1\leq i\leq n-1$. Let $\bbQ$ be an
ordered $r$-tuple which is obtained by fixing an order on $\bQ$.
Then for each $\ulam\in\mathcal{P}_n$, there are
$\H_{r,n}^{K}(q,a\!\bQ)$-module isomorphisms
$$\Bigl(\ts_{\bbQ}^{\lam}\Bigr)^{\sigma_a}\cong\ts_{\overrightarrow{a\!\bQ}}^{\lam},\,\,\,\,
\Bigl(\td_{\bbQ}^{\lam}\Bigr)^{\sigma_a}\cong\td_{\overrightarrow{a\!\bQ}}^{\lam},$$
where $\overrightarrow{a\!\bQ}$ denotes the ordered $r$-tuple which
is obtained from $\bbQ$ by multiplying $a$ on each component. In
particular, $\td_{\bbQ}^{\ulam}\neq 0$ if and only if
$\td_{\overrightarrow{a\!\bQ}}^{\ulam}\neq 0$. Moreover, for each
$\ulam\in\mathcal{K}_n$, we have
$\rm{h}_{\bbQ}(\ulam)=\rm{h}_{a\!\bbQ}(\ulam)$.
\end{lem}

\begin{proof} This follows directly from the definition of $\ts_{\bbQ}^{\ulam}$ and
$\td_{\bbQ}^{\ulam}$.
\end{proof}

\begin{lem} {\rm (\cite[Lemma 3.5]{Hu5})} \label{NT} Let K be a field which contains a primitive p-th root of unity $\eps$.
Suppose $p=dk$, where $p, d, k\in\mathbb{N}$, $\xi\in K$ is a primitive $d$-th root of unity. Then
there exists a primitive $p$-th root of unity $\zeta\in K$ such that $\zeta^k=\xi$.
\end{lem}

Using the above two lemmas, we can divide the proof of Theorem
\ref{main2} into the following two cases:
\smallskip

\noindent {\it Case 1.}\,
$q^{\mathbb{Z}}\cap\eps^{\mathbb{Z}}=\{1\}$,
$\bbQ=(\bbQ^{[1]},\cdots,\bbQ^{[p]})$, where for each $1\leq j\leq
p$, $$
\bbQ^{[j]}=\bigl(\eps^{j-1}q^{v_1},\cdots,\eps^{j-1}q^{v_{d}}\bigr),
$$ for some integers $0\leq v_1,\cdots,v_d\leq e-1$.

In this case, for each $1\leq i\leq p$, let $$
\lam^{[i]}:=\Bigl(\lam^{(d(i-1)+1)},\cdots,\lam^{(di)}\Bigr),\quad n_i:=|\lam^{[i]}|.
$$
Applying, in turn, \cite[Lemma 2.4, Corollary 5.6]{Hu5} and
\cite[(5.12)]{V}, we have that $$\begin{aligned}
&\quad\,\,\Bigl(\td_{\bbQ}^{\lam}\Bigr)^{\sigma}\cong\td_{(\eps\!\bbQ^{[1]},\eps\!\bbQ^{[2]},\cdots,
\eps\!\bbQ^{[p]})}^{\lam}\cong\td_{(\bbQ^{[2]},\cdots,\bbQ^{[p-1]},\bbQ^{[1]})}^{\lam}\\
&\cong\Ind_{\H_{n_1}^{\aff}\otimes\cdots\otimes\H_{n_{p-1}}^{\aff}\otimes\H_{n_{p}}^{\aff}}^{\H_n^{\aff}}\Bigl(
D^{\lam^{[1]}}_{\bbQ^{[2]}}\otimes\cdots\otimes
D^{\lam^{[p-1]}}_{\bbQ^{[p]}}\otimes D^{\lam^{[p]}}_{\bbQ^{[1]}}\Bigr)\\
&\cong\Ind_{\H_{n_p}^{\aff}\otimes\H_{n_{1}}^{\aff}\otimes\cdots\otimes\H_{n_{p-1}}^{\aff}}^{\H_n^{\aff}}\Bigl(
D^{\lam^{[p]}}_{\bbQ^{[1]}}\otimes D^{\lam^{[1]}}_{\bbQ^{[2]}}\otimes\cdots\otimes
D^{\lam^{[p-1]}}_{\bbQ^{[p]}}\Bigr)\\
&\cong \td_{\bbQ}^{(\lam^{[p]},\lam^{[1]},\cdots,
\lam^{[p-1]})},\end{aligned}
$$
It follows that
\addtocounter{equation}{6}\begin{equation}\label{equa45}
\HH_{\bbQ}(\lam)=\bigl(\lam^{[p]},\lam^{[1]},\cdots,\lam^{[p-1]}\bigr).
\end{equation}
Note that $\bbQ^{[1]},\cdots,\bbQ^{[p]}$ are $p$ different $q$-orbits. Therefore, in this case Theorem
\ref{main2} follows easily from (\ref{equa45}).
\medskip

\noindent {\it Case 2.}\, $\bbQ=(\bbQ^{[1]},\cdots,\bbQ^{[k]})$,
where $p={d}_0k$, $q$ is a primitive ${d}_0\ell$-th root of unity,
$q^{\ell}=\eps^{k}$ is a primitive ${d}_0$-th root of unity, and
$1\leq k<p$ is the smallest positive integer such that $\eps^{k}\in
q^{\mathbb{Z}}$, and for each $1\leq j\leq k$,
$$\begin{aligned}
\bbQ^{[j]}&=\bigl(\eps^{j-1}q^{v_1},\cdots,\eps^{j-1}q^{v_d},\eps^{k+j-1}q^{v_1},\cdots,
\eps^{k+j-1}q^{v_d},\cdots,\cdots\\
&\qquad\qquad\eps^{({d}_0-1)k+j-1}q^{v_1},\cdots,\eps^{({d}_0-1)k+j-1}q^{v_{d}}\bigr),
\end{aligned}
$$
where $0\leq v_1\leq v_2\leq\cdots\leq v_d<\l$ are some integers independent
of $j$. In particular, in each $q$-orbit $\bbQ^{[j]}$, we have (compare this with our assumption in the third paragraph in Section 3)
$$\begin{aligned}
&0\leq v_1\leq v_2\leq\cdots\leq v_d<\l+v_1\leq \l+v_2\leq\cdots\leq \l+v_d<\cdots\\
&\qquad <(d_0-1)\l+v_1\leq (d_0-1)\l+v_2\leq\cdots\leq (d_0-1)\l+v_d<e.\end{aligned}
$$

Note that, in this case by assumption $e=d_0l<\infty$. We actually have
two different approaches. The first one is based on the same
arguments in \cite[Section 4]{Hu5}, we leave the details to the
interested readers. In this paper, we adopt a second approach, which
is based on \cite[Proposition 2.10]{GJ}) and the following two
results. Note also that we have fixed an order of the parameters in $\bQ$. This
order is important in the following lemma.

\addtocounter{lem}{1}
\begin{lem} \label{main11} We keep the same assumption as in {\it Case 2}, and take $z_i=\eps^{i-1}$ for each integer $i$ with
$1\leq i\leq k=s$. Then the notion of FLOTW $r$-partition with respect to $(q,\bbQ)$
is well-defined. For any FLOTW $r$-partition $\lam$ with respect to $(q,\bbQ)$, and
any two nodes $\gamma=(a,b,c)$, $\gamma'=(a',b',c')$ of $\lam$ with
the same residue,  $$ b-a>b'-a'\iff \text{$\gamma$ is below
$\gamma'$ with respect to the FLOTW order.}
$$
\end{lem}

\begin{proof} Let $\gamma=(a,b,c), \gamma'=(a',b',c')$ be two nodes of $\lam$ with
the same residue. Then $Q_cq^{b-a}=Q_{c'}q^{b'-a'}$. By our assumption on the $r$-tuple $\bbQ$. We have that
$$
Q_c=\eps^{j-1+c_1k}q^{v_{c_2}},\quad Q_{c'}=\eps^{j-1+c'_1k}q^{v_{c'_2}},
$$
for some integers $1\leq j\leq k$, $0\leq c_1, c'_1\leq d_0-1$, $1\leq c_2, c'_2\leq d$. Therefore, $$
b-a+c_1\l+v_{c_2}=b'-a'+c'_1\l+v_{c'_2}+le
$$
for some integer $l$.

Suppose $b-a>b'-a'$. By assumption, $$
|(c_1\l+v_{c_2})-(c'_1\l+v_{c'_2})|<(d_0-1)\l+\l=e.
$$ It follows that we must have $l\geq 0$. If $l>0$, then $$
b-a+c_1\l+v_{c_2}>b'-a'+c'_1\l+v_{c'_2},
$$
which implies that $\gamma$ is below $\gamma'$ with respect to the FLOTW order;
while if $l=0$, then we must have $$
0<(b-a)-(b'-a')<e,\quad -e<(c_1\l+v_{c_2})-(c'_1\l+v_{c'_2})<0.
$$
Note that $|v_{c_2}-v_{c'_2}|<\l$. It follows that either $c_1<c'_1$ or $c_1=c'_1$ and $c_2<c'_2$. By the definition of our $\bbQ$, we deduce that $c<c'$.
Hence $\gamma$ is again below $\gamma'$ with respect to the FLOTW order. \smallskip

Conversely, if $\gamma$ is below $\gamma'$ with respect to the FLOTW order, it is also easy to deduce that $b-a>b'-a'$.
\end{proof}

Let $\omega$ be the permutation on $\{1,2,\cdots,r\}$ which is
defined by $$\begin{aligned}
xd_0d+y&\mapsto (x+1)d_0d+y,\quad \forall\,0\leq x<k-1,\, 1\leq y\leq d_0d,\\
(k-1)d_0d+\hat{x}&\mapsto d+\hat{x},\quad \forall\,1\leq
\hat{x}\leq d_0d-d,\\
r-d+\hat{y}&\mapsto \hat{y},\quad \forall\,1\leq \hat{y}\leq d.
\end{aligned}
$$
By the definition of our $\bbQ$ (see the second paragraph above
Lemma \ref{main11}), it is easy to see that
$\eps\!\bbQ=(Q_{\omega(1)},\cdots,Q_{\omega(r)})$. For any
$\mu=(\mu^{(1)},\cdots,\mu^{(r)})\in\mathcal{P}_n$, we define
$\omega(\mu):=(\mu^{(\omega^{-1}(1))},$
$\cdots,\mu^{(\omega^{-1}(r))})$. Recall that (see Lemma \ref{lm32})
the map $\kappa$ restricts to a bijection from the set
$\mathcal{K}_n$ of Kleshchev $r$-multipartitions of $n$ with respect
to $(q,Q_1,\cdots,Q_r)$ to the set $\mathcal{F}_n$ of FLOTW
$r$-partitions of $n$ with respect to $(q,Q_1,\cdots,Q_r)$.

\begin{lem}\label{main12} Assume that $e<\infty$. Let $\lam\in\mathcal{K}_n$ be a Kleshchev $r$ multipartition of
$n$ with respect to $(q,Q_1,\cdots,Q_r)$. Then, if
$\underline{\emptyset}\underset{r_1}{\twoheadrightarrow}\cdot
\underset{r_2}{\twoheadrightarrow}\cdot \cdots\cdots
\underset{r_n}{\twoheadrightarrow}\kappa(\lam) $ is a path from
$\underline{\emptyset}$ to $\kappa(\lam)$ in FLOTW's good lattice
with respect to $(q,Q_1,\cdots,$ $Q_r)$, then the sequence
$$\underline{\emptyset}\underset{\eps r_1 }{\twoheadrightarrow}\cdot
\underset{\eps r_2 }{\twoheadrightarrow}\cdot \cdots\cdots
\underset{\eps r_n }{\twoheadrightarrow}\omega(\kappa(\lam)) $$ also
defines a path in FLOTW's good lattice with respect to
$(q,Q_1,\cdots,Q_r)$, and it connects $\underline{\emptyset}$ to
$\omega(\kappa(\lam))$.
\end{lem}

\begin{proof} By \cite[Proposition 2.10]{GJ}, we know that for any FLOTW $r$-partition $\mu$,
$\omega(\mu)$ is also a FLOTW $r$-partition. Now let $\mu$ be a given FLOTW $r$-partition.
Lemma \ref{main11} implies that $\gamma=(a,b,c)$ is a good $x$-node of $\mu$ if and only if
$\omega(\gamma):=(a,b,\omega(c))$ is a good $\eps x$-node of $\omega(\mu)$, from which the lemma follows
immediately.
\end{proof}

We remark that, in the special case where $r=p$ and $\eps=q^l$,
Lemma \ref{main11} and \ref{main12} are proved in
\cite[(4.3.A)]{Ja0}.

\medskip
\noindent {\bf Proof of Theorem \ref{main2} in Case 2}:\, By
\cite[Proposition 2.10]{GJ} (see also the second line in Page 16 of \cite{GJ}), we know that for any Kleshchev $r$-multipartition
$\lam$,
$\kappa\bigl(\HH(\lam)\bigr)=\omega\bigl(\kappa(\lam)\bigr)$. Now
Theorem \ref{main2} in Case 2 follows from Lemma \ref{lm32} and
Lemma \ref{main12}. \qed
\medskip

We have the following result (compare \cite[Theorem 3.8]{Hu5}).
\addtocounter{prop}{9}
\begin{prop} \label{prop48} Assume that $e<\infty$, $p=d_0k$, $q$ is a primitive $d_0\l$-th root of unity, $q^{\l}=\eps^k$ is a primitive
$d_0$-th root of unity, and $1\leq k<p$ is the smallest
positive integer such that $\eps^k\in q^{\mathbb{Z}}$. For each
$1\leq j\leq k$, let
$$\begin{aligned}
\bbQ^{[j]}&=\bigl(\eps^{j-1}q^{v_1},\cdots,\eps^{j-1}q^{v_{d}},\eps^{k+j-1}q^{v_1},\cdots,\eps^{k+j-1}q^{v_{d}},\cdots,
\cdots,\\
&\qquad\qquad \eps^{({d}_0-1)k+j-1}q^{v_1},\cdots,\eps^{({d}_0-1)k+j-1}q^{v_{d}}\bigr),\end{aligned}
$$
Let $\bbQ:=(\bbQ^{[1]},\cdots,\bbQ^{[k]})$, the concatenation of ordered tuples. Let
$\lam=(\lam^{[1]},\cdots,\lam^{[k]})\in\mathcal{K}_n$, where for each $1\leq i\leq k$, $$
\lam^{[i]}=(\lam^{((i-1)d_0d+1)},\cdots,\lam^{((i-1)d_0d+d_0d)}).
$$ Then $$
\HH(\lam)=\Bigl(\HH'(\lam^{[k]}),\lam^{[1]},\cdots,\lam^{[k-1]}\Bigr),
$$
where $\HH'$ is an automorphism on the set of Kleshchev $d_0d$-multipartition of $n_k$ with respect to $(q,\bbQ^{\vee})$,
where $$
\bbQ^{\vee}:=\Bigl(q^{v_1},\cdots,q^{v_d},\eps^kq^{v_1},\cdots, \eps^kq^{v_d},\cdots,\cdots,
\eps^{(d_0-1)k}q^{v_1},\cdots,\eps^{(d_0-1)k}q^{v_d}\Bigr),$$ $\HH'$ is defined (similar to $\HH$) in the context of the cyclotomic Hecke algebra $\H_{q, \bbQ^{\vee}}(d_0d,d_0,n_k)$, where $n_k=|\lam^{[k]}|$.
\end{prop}

\begin{proof} Note that $\bbQ^{[1]},\cdots,\bbQ^{[k]}$ are $k$
different $q$-orbits. By \cite[(4.11)]{DM}, the assumption that
$\lam\in\mathcal{K}_n$ implies that for each $1\leq i\leq k$,
$\lam^{[i]}$ is a Kleshchev $d_0d$-multipartition with respect to
$(q,\bbQ^{[i]})$. In particular, this in turn implies that $$
\Bigl(\HH'(\lam^{[k]}),\lam^{[1]},\cdots,\lam^{[k-1]}\Bigr)\in\mathcal{K}_n.
$$
By the same reasoning, it is easy to see that we can find a path of
the following form $$\begin{aligned}
\underline{\emptyset}&=(\underbrace{\emptyset,\cdots,\emptyset}_{\text{$r$
copies}})\overset{r_1}{\twoheadrightarrow}\cdot\cdots
\overset{r_{n_1}}{\twoheadrightarrow}(\lam^{[1]},\underbrace{\emptyset,\cdots,\emptyset}_{\text{$r-d_0d$
copies}})\overset{r_{n_1+1}}{\twoheadrightarrow}\cdot\cdots\\
&\quad
\overset{r_{n_1+n_2}}{\twoheadrightarrow}(\lam^{[1]},\lam^{[2]},\underbrace{\emptyset,\cdots,\emptyset}_{\text{$r-2d_0d$
copies}})\overset{r_{n_1+n_2+1}}{\twoheadrightarrow}\cdot\cdots\overset{r_{n}}{\twoheadrightarrow}\lam
\end{aligned}$$ in Kleshchev's good lattice with respect
to $(q,\bbQ)$ which connects $\underline{\emptyset}$ to $\lam$.
Now we apply Theorem \ref{main2}, the proposition follows at once.
\end{proof}

\bigskip\bigskip
\section{Explicit formulas for the number of simple modules}
\smallskip\smallskip

In this section, we shall derive explicit formulas for the number
of simple modules over the cyclotomic Hecke algebras of type
$G(r,p,n)$.
\smallskip

By Lemma \ref{lm24} and \cite[(6.2), (6.3)]{Hu5},
\begin{equation}\label{equa62}\begin{split}
\#\Irr\bigl(\H_K(r,p,n)\bigr)&=\frac{1}{p}\Bigl\{\#\Irr\bigl(\H_{K}(r,n)\bigr)
-\sum_{1\leq m<p,m|p}N(m)\Bigr\}\\
&\qquad+\sum_{1\leq
m<p,m|p}\frac{N(m)}{m}\frac{p}{m},\end{split}\end{equation}
where for each integer $1\leq m,\widetilde{m}\leq p$ with $m|p, \widetilde{m}|p$, $$\begin{aligned}
N(\widetilde{m})&=\sum_{1\leq m\leq \widetilde{m},
m|\widetilde{m}}\mu(\widetilde{m}/m)\widetilde{N}(m)\\
\widetilde{N}(m):&=\#\bigl\{\lam\in\mathcal{K}_n\bigm|\HH^{m}(\lam)=\lam\bigr\},\end{aligned}
$$
and $\mu(?)$ is the  M\"{o}bius function. Note that (by \cite{AM}) $\#\Irr\bigl(\H_{K}(r,n)\bigr)$ is explicitly known.
Therefore, it suffices to compute $\widetilde{N}(m)$. By the discussion in last section, to compute $\widetilde{N}(m)$, it suffices to consider the following two cases:
\medskip

\noindent {\it Case 1.}\, $\bbQ=(\bbQ^{[1]},\cdots,\bbQ^{[p]})$,
$q^{\mathbb{Z}}\cap\eps^{\mathbb{Z}}=\{1\}$, and for each $1\leq
j\leq p$, $$
\bbQ^{[j]}=\bigl(\eps^{j-1}q^{v_1},\cdots,\eps^{j-1}q^{v_{d}}\bigr),
$$
where $0\leq v_1,\cdots,v_d\leq e-1$ are some integers independent
of $j$.
\medskip

\noindent {\it Case 2.}\, $\bbQ=(\bbQ^{[1]},\cdots,\bbQ^{[k]})$,
where $p={d}_0k$, $q$ is a primitive ${d}_0\ell$-th root of unity,
$e=d_0\l$, $q^{\ell}=\eps^{k}$ is a primitive ${d}_0$-th root of
unity, and $1\leq k<p$ is the smallest positive integer such that
$\eps^{k}\in q^{\mathbb{Z}}$, and for each $1\leq j\leq k$,
$$\begin{aligned}
\bbQ^{[j]}&=\bigl(\eps^{j-1}q^{v_1},\cdots,\eps^{j-1}q^{v_{d}},\eps^{k+j-1}q^{v_1},\cdots,\eps^{k+j-1}q^{v_{d}},\cdots,
\cdots,\\
&\qquad\qquad \eps^{({d}_0-1)k+j-1}q^{v_1},\cdots,\eps^{({d}_0-1)k+j-1}q^{v_{d}}\bigr),\end{aligned}
$$
where $0\leq v_1\leq\cdots\leq v_d\leq \l-1$ are some integers independent
of $j$.
\medskip

The following are the second two main results in this paper, which yield explicit
formulas for the number of simple modules.

\addtocounter{thm}{1}
\begin{thm} \label{mainthm3} With the notations and assumptions as in Case 1, let $1\leq m\leq p$ be an integer such that $m|p$. Let $\bbQ^{\vee}=\bigl(q^{v_1},\cdots,q^{v_{d}}\bigr)$. If $p\nmid mn$, then $\widetilde{N}(m)=0$; if\, $p|mn$, then $$
\widetilde{N}(m)=\sum_{n_1+\cdots+n_{m}=\frac{nm}{p}}\prod_{i=1}^{m}
\biggl(\#\Irr\H_{q,\bbQ^{\vee}}\Bigl(d,n_i\Bigr)
\biggr).$$
where $\H_{q,\bbQ^{\vee}}(d,n_i)$ is the Ariki--Koike
algebra with parameters $(q,\bbQ^{\vee})$ and of size $n_i$.
\end{thm}

\begin{thm} \label{mainthm4} With the notations and assumption as in Case 2, let $1\leq m\leq p$ be an integer such that
$m|p$. Let $a=\gcd(m,k)$, $\widetilde{d}=\gcd(\frac{m}{a},{d}_0)$.
Let $q''$ be a primitive $\widetilde{d}\l$-th root of unity. Let
$$\begin{aligned}
\bbQ^{\vee}&=\bigl((q'')^{v_1},\cdots,(q'')^{v_{d}},(q'')^{v_1+\l},\cdots,(q'')^{v_{d}+\l},\cdots,\cdots,\\
&\qquad\qquad (q'')^{v_1+(\widetilde{d}-1)\l},\cdots,(q'')^{v_{d}+(\widetilde{d}-1)\l}\bigr).\end{aligned}
$$
If $k\nmid na$, then $\widetilde{N}(m)=0$; if\, $k|na$, then $$
\widetilde{N}(m)=\sum_{n_1+\cdots+n_{a}=\frac{na}{k}}\prod_{i=1}^{a}
\biggl(\#\Irr\H_{q'',\bbQ^{\vee}}\Bigl(\widetilde{d}d,\frac{\widetilde{d}n_i}{{d}_0}\Bigr)
\biggr).$$
where $\H_{q'',\bbQ^{\vee}}(\widetilde{d}d,\frac{\widetilde{d}n_i}{{d}_0})$ is the Ariki--Koike
algebra with parameters $(q'',\bbQ^{\vee})$ and of size $\frac{\widetilde{d}n_i}{d_0}$, and the number
$\#\Irr\H_{q'',\bbQ^{\vee}}\Bigl(\widetilde{d}d,\frac{\widetilde{d}n_i}{{d}_0}\Bigr)$
is understood as $0$ if $d_0\nmid\widetilde{d}n_i$.
\end{thm}
\medskip

\noindent
{\bf Proof of Theorem \ref{mainthm3}}:\, Let $\lam=(\lam^{(1)},\cdots,\lam^{(r)})$. We write $\lam=(\lam^{[1]},\cdots,$ $\lam^{[p]})$, where for each $1\leq j\leq p$, $$
\lam^{[j]}=(\lam^{((j-1)d+1)},\cdots,\lam^{(jd)}).
$$
By \cite[(4.11)]{DM}, $\lam\in\mathcal{K}_n$ if and only for each
$1\leq i\leq p$, $\lam^{[i]}$ is a Kleshchev $d$-multipartition
with respect to $(q,\bbQ^{\vee})$.

By (\ref{equa45}), we know that $$
\HH^{m}(\lam)=\Bigl(\underbrace{\lam^{[p-m+1]},\lam^{[p-m+2]},\cdots,\lam^{[p]}}_{\text{$m$ terms}},\underbrace{\lam^{[1]},\lam^{[2]},\cdots,\lam^{[p-m]}}_{\text{$p-m$ terms}}\Bigr).
$$
It is easy to see that $\HH^m(\lam)=\lam$ if and only if $$
\lam^{[i]}=\lam^{[i+lm]}\quad\text{for each $1\leq l\leq p/m-1$ and each $1\leq i\leq m$,}
$$
from which Theorem \ref{mainthm3} follows immediately.
\bigskip

We now turn to the proof of Theorem \ref{mainthm4}. We will use the same strategy as in \cite[Section 5]{Hu5}, where
the proof makes use of Naito--Sagaki's work (\cite{NS2}, \cite{NS1}).

We keep the notations and assumption as in Case 2. For the moment, we assume that $k=1$. In other words, $\eps=q^{\ell}$,
$q$ is a primitive $p\ell$-th root of unity, $e=p\l$ and $$\begin{aligned}
\bbQ&=\bigl(q^{v_1},\cdots,q^{v_{d}},\eps q^{v_1},\cdots,\eps q^{v_{d}},\cdots,\cdots,\\
&\qquad\qquad\qquad \eps^{p-1}q^{v_1},\cdots,\eps^{p-1}q^{v_{d}}\bigr),\end{aligned}
$$
for some integers $0\leq v_1\leq\cdots\leq v_d\leq\l-1$. We consider the
affine Kac--Moody algebra $\mathfrak{g}=\ksl_{p\l}$ of type
$A_{p\l-1}^{(1)}$. Let $\mathfrak{h}$ be the Cartan subalgebra of $\mathfrak{g}$, let $W$ be the Weyl group of $\mathfrak{g}$. Let $I:=\mathbb{Z}/p\l\mathbb{Z}$. Let
$\pi:\,I\rightarrow I$ be the Dynkin diagram automorphism of order
$p/m$ defined by $\bar{i}=i+p\l\Z\mapsto
\bar{i}-\overline{m\l}=i-m\l+p\l\Z$ for any $\bar{i}\in I$. By \cite{FSS}, $\pi$ induces a Lie algebra automorphism (which is called the
diagram outer automorphism) $\pi\in\Aut(\mathfrak{g})$ of order $p/m$ and a linear automorphism
$\pi^{\ast}\in GL(\mathfrak{h}^{\ast})$ of order $p/m$. Let
$\check{\mathfrak{g}}$ be the corresponding orbit Lie algebra. Then
(by \cite[(6.4)]{Hu5})
$$
\check{\mathfrak{g}}=\begin{cases} \ksl_{m\l}, &\text{if $m\l>1$,}\\
\mathbb{C}, &\text{if $m=\l=1$.}
\end{cases}
$$
Let $\check{\mathfrak{h}}$, $\check{W}$,
$\{\check{\Lambda}_i\}_{0\leq i\leq m\l-1}$ denote the Cartan
subalgebra, the Weyl group, the set of fundamental dominant weights
of $\check{\mathfrak{g}}$ respectively. Let $\widetilde{W}=\{w\in W|\pi^{\ast}w=w\pi^{\ast}\}$.
There exists a linear
automorphism
$P_{\pi}^{\ast}:\,\check{\mathfrak{h}}^{\ast}\rightarrow
\bigl(\mathfrak{h}^{\ast}\bigr)^{\circ}:=\bigl\{\Lambda\in\mathfrak{h}^{\ast}
\bigm|\pi^{\ast}(\Lambda)=\Lambda\bigr\}$ and a group isomorphism
$\Theta:\,\check{W}\rightarrow\widetilde{W}$ such that
$\Theta(\check{w})=P_{\pi}^{\ast}\check{w}\bigl(P_{\pi}^{\ast}\bigr)^{-1}$
for each $\check{w}\in\check{W}$. By \cite[\S6.5]{FSS}, for each $0\leq
i<m\l$, $$
P_{\pi}^{\ast}(\check{\Lambda}_i)=\Lambda_{i}+\Lambda_{i+m\l}+\Lambda_{i+2m\l}+\cdots+\Lambda_{i+(p-m)\l}+C\delta,
$$
where $C\in\mathbb{Q}$ is some constant depending on $\pi$, $\delta$
denotes the null root of $\mathfrak{g}$. Let
$$\check{\Lambda}:=\sum_{i=1}^{d}\sum_{j=1}^{m}\Lambda_{v_i+(j-1)\l},\quad
\Lambda:=\sum_{i=1}^{d}\sum_{j=1}^{p}\Lambda_{v_i+(j-1)\l}.$$ Then
it follows that $P_{\pi}^{\ast}(\check{\Lambda})=\Lambda+C'\delta$,
for some $C'\in\mathbb{Q}$.

Let $\eps':=\eps^{p/m}$, which is a primitive $m$-th root of unity. By Lemma \ref{NT}, we can find a primitive $m\l$th root of unity $q'$, such that
$(q')^{\l}={\eps'}$. Let
$$\begin{aligned}
\bbQ^{\vee}&=\bigl(q^{v_1},\cdots,q^{v_{d}},\eps q^{v_1},\cdots, \eps q^{v_{d}},\cdots,\cdots,\\
&\qquad\qquad\qquad \eps^{m-1}q^{v_1},\cdots,\eps^{m-1}q^{v_{d}}\bigr),\end{aligned}
$$
By the same argument as in \cite[(6.9),(6.10)]{Hu5}, we get that

\addtocounter{cor}{3}
\begin{cor} \label{maincor} With the notation as above, there exists a bijection
$\eta:\,\check{\lam}\mapsto\lam$ from the set of Kleshchev
$dm$-multipartitions $\check{\lam}$ of $nm/p$ with respect to
$(q',\bbQ^{\vee})$ onto the set of
Kleshchev $dp$-multipartitions $\lam$ of $n$ with respect to
$(q,\bbQ)$ satisfying $\HH^m(\lam)=\lam$, such that if
$$
\underbrace{(\emptyset,\cdots,\emptyset)}_{\text{$dm$ copies
}}\overset{r_1}{\twoheadrightarrow}\cdot
\overset{r_2}{\twoheadrightarrow}\cdot \cdots\cdots
\overset{r_s}{\twoheadrightarrow}\check\lam $$ is a path from
$\underbrace{(\emptyset,\cdots,\emptyset)}_{\text{$dm$ copies}}$
to $\check\lam$ in Kleshchev's good lattice with respect to
$(q',\bbQ^{\vee})$, where $s:=nm/p$, then the sequence
$$\begin{aligned}&\underbrace{(\emptyset,\cdots,\emptyset)}_{\text{$dp$ copies}}
\overset{r_1}{\twoheadrightarrow}\cdot\overset{m\l+r_1}{\twoheadrightarrow}\cdot
\cdots\overset{(p-m)\l+r_1}{\twoheadrightarrow}\cdot
\overset{r_2}{\twoheadrightarrow}\cdot
\overset{m\l+r_2}{\twoheadrightarrow}\cdot
\cdots\overset{(p-m)\l+r_2}{\twoheadrightarrow}\cdot\\
&\qquad \qquad\cdots\,\cdot
\overset{r_s}{\twoheadrightarrow}\cdot\overset{m\l+r_s}{\twoheadrightarrow}\cdots
\overset{(p-m)\l+r_s}{\twoheadrightarrow}\lam\end{aligned}
$$ defines a path in Kleshchev's good lattice (w.r.t.,
$(q,\bbQ)$) satisfying $\HH^m(\lam)=\lam$. In particular, Theorem \ref{mainthm4} is valid in the case $k=1$. That is, $$
\widetilde{N}(m)=\#\Irr\H_{q',\bbQ^{\vee}}\Bigl(md,\frac{mn}{p}\Bigr).
$$
\end{cor}

Now we consider the case where $k>1$. We keep the notation and
assumption as in Case 2. That is,
$\bbQ=(\bbQ^{[1]},\cdots,\bbQ^{[k]})$, where $p={d}_0k$, $q$ is a
primitive ${d}_0\ell$-th root of unity, $q^{\ell}=\eps^{k}$ is a
primitive ${d}_0$-th root of unity, $e=d_0\l$, and $1<k<p$
is the smallest positive integer such that $\eps^{k}\in
q^{\mathbb{Z}}$, and for each $1\leq j\leq k$, $$\begin{aligned}
\bbQ^{[j]}&=\bigl(\eps^{j-1}q^{v_1},\cdots,\eps^{j-1}q^{v_{d}},\eps^{k+j-1}q^{v_1},\cdots,\eps^{k+j-1}q^{v_{d}},
\cdots,\cdots,\\
&\qquad\qquad \eps^{({d}_0-1)k+j-1}q^{v_1},\cdots,\eps^{({d}_0-1)k+j-1}q^{v_{d}}\bigr),\end{aligned}
$$
for some integers $0\leq v_1\leq\cdots\leq v_d\leq\l-1$. Let
$\lam=(\lam^{[1]},\cdots,\lam^{[k]})$, where for each $1\leq i\leq
k$, $$
\lam^{[i]}=(\lam^{((i-1)d_0d+1)},\cdots,\lam^{((i-1)d_0d+d_0d)}).
$$
Let $n_i:=|\lam^{[i]}|$ for each $1\leq i\leq k$. Clearly, $\lam\in\mathcal{K}_n$ if and only if for each $1\leq i\leq k$, $\lam^{[i]}\in\mathcal{K}_{n_i}$, where
$\mathcal{K}_{n_i}$ denotes the set of Kleshchev $d_0d$-multipartitions of $n_i$ with respect to
$$(q,\underbrace{q^{v_1},\cdots, q^{v_d}}_{\text{$d$ terms}},\underbrace{\eps^k q^{v_1},\cdots, \eps^k q^{v_d}}_{\text{$d$ terms}},\cdots,\underbrace{\eps^{({d}_0-1)k}q^{v_1},\cdots,
\eps^{({d}_0-1)k}q^{v_d}}_{\text{$d$ terms}}). $$

Suppose that $1\leq a\leq\min\{m,k\}$ is the greatest common divisor of $m$ and $k$. We define $$\begin{aligned}
&\widetilde{\Sigma}(k,m):=\biggl\{\bigl(\lam^{[1]},\cdots,\lam^{[a]}\bigr)\vdash
\frac{na}{k}\biggm|\begin{matrix}
\text{$\lam^{[i]}\in\mathcal{K}_{n_i}, (\HH')^{m/a}({\lam^{[i]}})={\lam^{[i]}}$,}\\
\text{$\forall\,1\leq i\leq a,\,\,\sum_{i=1}^{a}n_i=\frac{na}{k}$}\end{matrix}\biggr\},\\
&\widetilde{N}(k,m):=\#\widetilde{\Sigma}(k,m),\\
&\widetilde{\Sigma}(m):=\bigl\{\lam\in\mathcal{K}_n\bigm|\HH^{m}(\lam)=\lam\bigr\}.\end{aligned}
$$
where $\HH'$ is the same as in Proposition \ref{prop48}.

With the Proposition \ref{prop48} in mind, it is easy to see that
the same argument in the proof of \cite[Lemma 6.16]{Hu5} proves the
following result.

\addtocounter{lem}{4}
\begin{lem} \label{lm55} The map which sends
$\lam=\bigl(\lam^{[1]},\cdots,$ $\lam^{[k]}\bigr)$ to
$\overline{\lam}:=\bigl(\lam^{[1]},\cdots,$ $\lam^{[a]}\bigr)$
defines a bijection from the set $\widetilde{\Sigma}(m)$ onto the
set $\widetilde{\Sigma}(k,m)$.
\end{lem}

Let $\widetilde{d}:=\gcd(d_0,\frac{m}{a})$. Note that
$(\HH')^{d_0}(\lam^{[i]})=\lam^{[i]}$ for each $1\leq i\leq k$.
Therefore, \addtocounter{equation}{4}\begin{equation}\label{equa56}
\text{$(\HH')^{m/a}(\lam^{[i]})=\lam^{[i]}$ if and only
$(\HH')^{\widetilde{d}}(\lam^{[i]})=\lam^{[i]}$.}
\end{equation}
Note that we have just proved Theorem \ref{mainthm4} in the case
where $k=1$, it is now easy to see that Theorem \ref{mainthm4} in
the case where $k>1$ follows directly from Proposition
\ref{prop48}, Corollary \ref{maincor}, Lemma \ref{lm55} and (\ref{equa56}). This completes
the proof of Theorem \ref{mainthm4} in all cases.
\bigskip\bigskip

\bigskip\bigskip
\bigskip\bigskip

\bibliographystyle{amsplain}

\bigskip\bigskip
\bigskip\bigskip\bigskip\bigskip

\centerline{\bf Acknowledgement}\bigskip

This work was supported by National Natural Science Foundation of
China (Project 10401005) and by the Program NCET and by the URF of
Victoria University of Wellington.
\bigskip\bigskip\bigskip

\vfill
\noindent
Department of Applied Mathematics,\\[2pt]
\noindent
Beijing Institute of Technology\\[2pt]
\noindent
Beijing, 100081\\[2pt]
The People's Republic of China\\[3pt]
\noindent
E-mail: junhu303@yahoo.com.cn
\bigskip\bigskip

\eject

\section*{\bf Appendix by Xiaoyi Cui}
\bigskip

Let $R$ be a commutative ring with identity $1_{R}$. Let $B$ be a
finitely generated $R$-free $R$-algebra. Let $A$ be a $\Z/rZ$-graded
algebra over the subalgebra $B$ with grading
$A=\oplus_{j=0}^{r-1}a^jB$ for a unit $a\in A$, $aB=Ba$, $a^r\in B$.
Furthermore, we assume that $r\cdot 1_{R}$ is a unit in $R$, and $R$
contains a primitive $r$th root of unity $\epsilon$.
\smallskip

Let $\bigl\{b_i\bigr\}_{i=1}^{s}$ be an $R$-basis of $B$. Then, the
set $\bigl\{b_ia^j\bigr\}_{1\leq i\leq s, 0\leq j<r}$ is an
$R$-basis of $A$. Furthermore, the set $$
\Bigl\{b_{i_1}a^{i_2}\otimes_{B}a^{i_3}\Bigm|1\leq i_1\leq s, 0\leq
i_2,i_3<r\Bigr\}
$$
is an $R$-basis of $A\otimes_{B}A$.\smallskip

Let $\sigma$ be the automorphism of $A$ which is defined by $$
a^jx\mapsto\epsilon^ja^j x,\,\,\forall\,x\in B,\, j\in\mathbb{Z}.
$$
For each integer $j$ with $0\leq j<r$, recall that $A^{\sigma^j}=A$
as $R$-module, the left $A$-action on $A^{\sigma^j}$ is given by the
usual left multiplication, while the right $A$-action on
$A^{\sigma^j}$ is given by the usual right multiplication twisted by
$\sigma^j$. To avoid confusion, we use the symbol
$\bigl(b_{i_1}a^{i_2}\bigr)_{(j)}$ to denote the element
$b_{i_1}a^{i_2}$ in $A^{\sigma^j}$. Then the elements
$$ \bigl(b_{i_1}a^{i_2}\bigr)_{(j)},\,\,\text{where
$i_1\in\{1,2,\cdots,s\}$, $i_2,j\in\{0,1,\cdots,r-1\}$,}
$$
form an $R$-basis of $\oplus_{j=0}^{r-1}A^{\sigma^j}$.
\medskip

The conclusion of the following fact is contained in the proof of
[G, Proposition 2.2]. However, the argument given by the
proof of [G, Proposition 2.2] contains a gap. That is, in the
10th line of Page 527, Genet's claim about the determinant of the
representing matrix is generally false. In fact, it is a quite
nontrivial job to calculate the the determinant of that representing
matrix, as one can see from the following proof.
\medskip\smallskip

\noindent
{\bf Fact:}\,\,  Let $\varphi:A\otimes_{B}A\rightarrow
\oplus_{j=0}^{r-1}A^{\sigma^j}$ be the $R$-linear homomorphism
defined on basis by $$
\varphi\Bigl(b_{i_1}a^{i_2}\otimes_{B}a^{i_3}\Bigr)=\oplus_{j=0}^{r-1}\Bigl(\epsilon^{ji_3}\bigl(b_{i_1}a^{i_2+i_3}
\bigr)_{(j)}\Bigr),
$$
where $i_1\in\{1,2,\cdots,s\}$, $i_2,i_3\in\{0,1,\cdots,r-1\}$. Then
$\varphi$ is an ($A,A$)-bimodule isomorphism.

\smallskip
\begin{proof} It is easy to verify that $\varphi$ is an ($A,A$)-bimodule
homomorphism. Therefore, it remains to show that $\varphi$ is an
$R$-linear isomorphism.\smallskip

For any integers $i_1,i_2,i_3$ with $1\leq i_1\leq s, 0\leq
i_2,i_3\leq r-1$, we define $$ X_{i_3 sr+i_2
s+i_1}:=b_{i_1}a^{i_2}\otimes_{B}a^{i_3}.
$$
Then the set $\bigl\{X_1,X_2,\cdots,X_{sr^2}\bigr\}$ is an ordered
$R$-basis of $A\otimes_{B}A$.\smallskip

For any integers $k_1,k_2,k_3$ with $1\leq k_2\leq s, 0\leq
k_1,k_3\leq r-1$, we define $$ Y_{k_1 sr+k_3
s+k_2}:=\bigl(b_{k_2}a^{k_3}\bigr)_{(k_1)}.
$$
Then the set $\bigl\{Y_1,Y_2,\cdots,Y_{sr^2}\bigr\}$ is an ordered
$R$-basis of $\oplus_{j=0}^{r-1}A^{\sigma^j}$.\smallskip

We want to compute the determinant of the representing matrix of
$\varphi$ with respect to the ordered $X$-basis and the ordered
$Y$-basis. Suppose that $$
\varphi\Bigl(X_1,X_2,\cdots,X_{sr^2}\Bigr)=\Bigl(Y_1,Y_2,\cdots,Y_{sr^2}\Bigr)M,
$$
where $M$ is an $sr^2\times sr^2$ matrix. For any integers $i,j$
with $0\leq i,j\leq r-1$, we use $X\!\!\downarrow_{i_3=i}$ (resp.,
$Y\!\!\downarrow_{k_1=j}$) to denote the naturally ordered basis
elements
$$ \bigl\{X_{i sr+i_2 s+i_1}\bigm|1\leq i_1\leq s, 0\leq i_2\leq
r-1\bigr\}
$$
$\Bigl($ resp., $\bigl\{Y_{jsr+k_3 s+k_2}\bigm|0\leq k_2,k_3\leq
r-1\bigr\}.\,\Bigr)$

We use $\iota_i$ to denote the embedding from the free $R$-submodule
spanned by elements in $X\!\!\downarrow_{i_3=i}$ into
$A\otimes_{B}A$; we use $p_j$ to denote the natural projection from
$\oplus_{t=0}^{r-1}A^{\sigma^t}$ onto $A^{\sigma^j}$ (i.e., the free
$R$-submodule spanned by the elements in
$Y\!\!\downarrow_{k_1=j}$.\smallskip

Let
$$\varphi\!\!\downarrow_{i_3=i, k_1=j}:=p_j\circ\varphi\circ\iota_i.
$$
Note that $X\!\!\downarrow_{i_3=i}$ (resp.,
$Y\!\!\downarrow_{k_1=j}$) is a consecutive part of
$(X_1,X_2,\cdots,X_{sr^2})$ (resp., of $(Y_1,Y_2,\cdots,Y_{sr^2})$).
Therefore, we can partition $M$ as
follows. $$M=\left(\begin{matrix} M_{0,0} & M_{0,1} & \cdots & M_{0,r-1}\\
M_{1,0} & M_{1,1} & \cdots & M_{1,r-1}\\
\vdots & \vdots & \vdots & \vdots\\
M_{r-1,0} & M_{r-1,1} & \cdots & M_{r-1,r-1}\\
\end{matrix}\right),
$$
where for each pair of integers $(i,j)$ with $0\leq i,j\leq r-1$,
$M_{j,i}$ is the representing matrix of
$\varphi\!\!\downarrow_{i_3=i, k_1=j}$ with respect to the ordered
basis $X\!\!\downarrow_{i_3=i}$ and $Y\!\!\downarrow_{k_1=j}$. That
is $$ \varphi\!\!\downarrow_{i_3=i,
k_1=j}\Bigl(X\!\!\downarrow_{i_3=i}\Bigr)=\Bigl(Y\!\!\downarrow_{k_1=j}\Bigr)M_{j,i}.
$$

Note that each $M_{j,i}$ is an $rs\times rs$ matrix. We claim
that\begin{enumerate} \item for each pair of integers $(i,j)$ with
$0\leq i,j\leq r-1$, $$ M_{j,i}=\epsilon^{ji}M_{0,i}.
$$
\item for each integer $i$ with $0\leq i\leq r-1$,
$M_{0,i}=(M_{0,1})^i$.
\end{enumerate}\smallskip

In fact, claim (1) follows directly from the definition of $\varphi$
and our ordering of the $X$ basis and the $Y$ basis. It suffices to
prove the claim (2).\smallskip

For any integer $i$ with $0\leq i\leq r-1$, we set $$
\varphi_i:=\varphi\!\!\downarrow_{i_3=i,
k_1=0}:\,\,b_{i_1}a^{i_2}\otimes a^i\mapsto
\bigl(b_{i_1}a^{i_2+i}\bigr)_{(0)}.
$$
We identify the free $R$-submodule of $A\otimes_{B}A$ spanned by
elements in $X\!\!\downarrow_{i_3=i}$ with $A^{\sigma^0}=A$ via $$
b_{i_1}a^{i_2}\otimes a^i\leftrightarrow
\bigl(b_{i_1}a^{i_2}\bigr)_{(0)},\,\,\,\text{for any integers $1\leq
i_1\leq s$, $0\leq i_2\leq r-1$.}
$$
With the above identification in mind, it is easy to see that $$
\varphi_i=\bigl(\varphi_1\bigr)^{i}.
$$
As a result, $M_{0,i}=(M_{0,1})^i$. This proves claim (2).\smallskip

Therefore, $$\begin{aligned} M&=\begin{pmatrix} I & M_{0,1} &
(M_{0,1})^2 & \cdots & (M_{0,1})^{r-1}\\
I & \epsilon M_{0,1} &
(\epsilon M_{0,1})^2 & \cdots & (\epsilon M_{0,1})^{r-1}\\
I & \epsilon^2 M_{0,1} &
(\epsilon^2 M_{0,1})^2 & \cdots & (\epsilon^2 M_{0,1})^{r-1}\\
\vdots & \vdots &\vdots &\vdots & \vdots \\
I & \epsilon^{r-1}M_{0,1} & (\epsilon^{r-1}M_{0,1})^2 & \cdots &
(\epsilon^{r-1}M_{0,1})^{r-1}
\end{pmatrix}\\
&=\begin{pmatrix} I & I &
(I)^2 & \cdots & (I)^{r-1}\\
I & \epsilon I &
(\epsilon I)^2 & \cdots & (\epsilon I)^{r-1}\\
I & \epsilon^2 I &
(\epsilon^2 I)^2 & \cdots & (\epsilon^2 I)^{r-1}\\
\vdots & \vdots &\vdots &\vdots & \vdots \\
I & \epsilon^{r-1}I & (\epsilon^{r-1}I)^2 & \cdots & (\epsilon^{r-1}I)^{r-1}
\end{pmatrix}\times\\
&\qquad\qquad
\begin{pmatrix} I & 0 &
0 & \cdots & 0\\
0 & M_{0,1} &
0 & \cdots & 0\\
0 & 0 &
M_{0,1}^2 & \cdots & 0\\
\vdots & \vdots &\vdots &\vdots & \vdots \\
0 & 0 & 0 & \cdots & (M_{0,1})^{r-1}
\end{pmatrix}\\
&=V_r\times D,
\end{aligned}
$$
where $I$ is the $rs\times rs$ identity matrix, $$\begin{aligned}
V_r&=\begin{pmatrix} I & I &
I & \cdots & I\\
I & \epsilon I &
\epsilon^2I & \cdots & \epsilon^{r-1}I\\
I & \epsilon^2 I &
\epsilon^4 I & \cdots & \epsilon^{2(r-1)}I\\
\vdots & \vdots &\vdots &\vdots & \vdots \\
I & \epsilon^{r-1}I & \epsilon^{2(r-1)}I & \cdots & \epsilon^{(r-1)^2}I
\end{pmatrix};\\
D&=\begin{pmatrix} I & 0 &
0 & \cdots & 0\\
0 & M_{0,1} &
0 & \cdots & 0\\
0 & 0 &
M_{0,1}^2 & \cdots & 0\\
\vdots & \vdots &\vdots &\vdots & \vdots \\
0 & 0 & 0 & \cdots & (M_{0,1})^{r-1}
\end{pmatrix}.
\end{aligned}
$$
Hence, $\det M=\det V_r\det D=\det V_r\bigl(\det
M_{0,1}\bigr)^{r(r-1)/2}$. To show that $\varphi$ is an isomorphism,
it suffices to show that $\det M$ is a unit in $R$. Therefore, it
suffices to show that both $\det V_r$ and $\det M_{0,1}$ are units
in $R$.\smallskip

By assumption, $a$ is invertible in $A$, which implies that the
elements in $\bigl\{b_1a^r,b_2a^r,\cdots,b_sa^r\bigr\}$ are
$R$-linear independent. Also, the condition that $a^r\in B$ implies
that there exists a matrix $C=(C_{i,j})_{s\times s}\in M_{s\times
s}(R)$, such that $\bigl(b_1a^r,b_2a^r,\cdots,b_sa^r\bigr)=(b_1,b_2,\cdots,b_s)C$.
\smallskip

By the condition that $a$ is invertible in $A$ we deduce that $C$ is
invertible in $M_{s\times s}(R)$. In particular, $\det C$ is
invertible in $R$, i.e., a unit in $R$. By direct calculation, we
know that $$ M_{0,1}=\begin{pmatrix} 0 & 0 &
0 & \cdots & 0 & C\\
I & 0 & 0 & \cdots & 0 & 0\\
0 & I & 0 & \cdots & 0 & 0\\
\vdots & \vdots &\vdots &\vdots & \vdots & \vdots \\
0 & 0 & 0 & \cdots & I & 0
\end{pmatrix},
$$
where each $0$ denotes an $s\times s$ zero matrix, each $I$ denote
an $s\times s$ identity matrix. As a consequence, $\det
M_{0,1}=(-1)^{(r-1)s^2}\det C$ is invertible in $R$. \smallskip

It remains to show that $\det V_r$ is invertible in $R$. In fact,
$$
\begin{aligned}
&\quad\,\det V_r\\
&=\left|
\begin{matrix} I & I &
I & \cdots & I\\
0 & (\epsilon-1)I &
(\epsilon^2-1)I & \cdots & (\epsilon^{r-1}-1)I\\
0 & (\epsilon^2-\epsilon)I &
(\epsilon^4-\epsilon^2)I & \cdots & (\epsilon^{2(r-1)}-\epsilon^{r-1})I\\
\vdots & \vdots &\vdots &\vdots & \vdots \\
0 & (\epsilon^{r-1}-\epsilon^{r-2})I & (\epsilon^{2(r-1)}-\epsilon^{2(r-2)})I &
\cdots & (\epsilon^{(r-1)^2}-\epsilon^{(r-1)(r-2)})I
\end{matrix}\right|.
\end{aligned}
$$
Hence $$
\begin{aligned}
&\quad\,\det V_r\\
&=\left|
\begin{matrix} (\epsilon-1)I &
(\epsilon^2-1)I & \cdots & (\epsilon^{r-1}-1)I\\
(\epsilon^2-\epsilon)I &
(\epsilon^4-\epsilon^2)I & \cdots & (\epsilon^{2(r-1)}-\epsilon^{r-1})I\\
\vdots &\vdots &\vdots & \vdots \\
(\epsilon^{r-1}-\epsilon^{r-2})I & (\epsilon^{2(r-1)}-\epsilon^{2(r-2)})I & \cdots &
(\epsilon^{(r-1)^2}-\epsilon^{(r-1)(r-2)})I
\end{matrix}\right|\\
&=\bigl((\epsilon-1)(\epsilon^2-1)\cdots(\epsilon^{r-1}-1)\bigr)^{rs}\left|
\begin{matrix} I &
I & \cdots & I\\
\epsilon I &
\epsilon^2I & \cdots & \epsilon^{r-1}I\\
\epsilon^2I &
\epsilon^4I & \cdots & \epsilon^{2(r-1)}I\\
\vdots &\vdots &\vdots & \vdots \\
\epsilon^{r-2}I & \epsilon^{2(r-2)}I & \cdots & \epsilon^{(r-1)(r-2)}I
\end{matrix}\right|\\
&=\Bigl(\epsilon^{(r-1)(r-2)/2}\prod_{t=1}^{r-1}(\epsilon^t-1)\Bigr)^{rs}\left|
\begin{matrix} I &
I & \cdots & I\\
I & \epsilon I & \cdots & \epsilon^{r-2}I\\
I & \epsilon^2I & \cdots & \epsilon^{2(r-2)}I\\
\vdots &\vdots &\vdots & \vdots \\
I & \epsilon^{r-2}I & \cdots & \epsilon^{(r-2)^2}I
\end{matrix}\right|\end{aligned}$$

$$
=\Bigl(\epsilon^{(r-1)(r-2)/2}\prod_{t=1}^{r-1}(\epsilon^t-1)\Bigr)^{rs}\det
V_{r-1}.
$$
Now by an easy induction argument, it is easy to see that $$ \det
V_r=\Bigl(\epsilon^{\sum_{a=1}^{r-2}\frac{a(a+1)}{2}}\prod_{b=1}^{r-1}\prod_{t=1}^{b}(\epsilon^t-1)
\Bigr)^{rs}.
$$
Note that $r=\prod_{1\leq j\leq r-1}(1-\epsilon^j)$. By assumption, $r$
is a unit in $R$. It follows that $\det V_r$ must be an invertible
element in $R$. This completes the proof.

\end{proof}

\bigskip\bigskip

\bigskip\bigskip\bigskip\bigskip
\bibliographystyle{amsplain}

\begin{thebibliography}{99}


\bibitem{A1} S. Ariki, On the decomposition numbers of the Hecke algebra of $G(m,1,n)$,
   {\em J. Math. Kyoto Univ.} {\bf 36} (1996), 789--808.

\bibitem{A2} S. Ariki, On the classification of simple modules
for cyclotomic Hecke algebras of type $G(m,1,n)$ and Kleshchev
multi-partitions, {\em Osaka J. Math.} (4) {\bf 38} (2001), 827--837.

\bibitem{A3} S. Ariki, Representation theory of a Hecke algebra of type $G(r,p,n)$,
   {\em J. Algebra} {\bf 177} (1995), 164--185.

\bibitem{A4} S. Ariki, Proof of the modular branching rule for cyclotomic Hecke algebras,
{\em J. Algebra} {\bf 306} (2006), 290--300.

\bibitem{AJ} S. Ariki and N. Jacon, Dipper--James--Murphy's conjecture for Hecke algebras of type B,
preprint, math.RT/0703447.

\bibitem{AK} S. Ariki and K. Koike, A Hecke algebra of
$(\Z/r\Z)\wr\BS_n$ and construction of its representations, {\em Adv.
Math.} {\bf 106} (1994), 216--243.

\bibitem{AM} S. Ariki and A. Mathas, The number of simple
modules of the Hecke algebras of type $G(r,1,n)$, {\em Math. Z.} (3)
{\bf 233} (2000), 601--623.

\bibitem{BM} M. Brou\'e and G. Malle, Zyklotomische
Heckealgebren, {\em Ast\'erisque} {\bf 212} (1993), 119--189.

\bibitem{DJ} R. Dipper and G. D. James, Representations of Hecke algebras of general linear groups,
{\em Proc. London. Math. Soc.} (3) {\bf 52} (1986), 20--52.

\bibitem{DJM} R. Dipper, G.D. James and A. Mathas, Cyclotomic
$q$-Schur algebras, {\em Math. Z.} {\bf 229} (1998), 385--416.

\bibitem{DM} R. Dipper and A. Mathas, Morita equivalence of
Ariki-Koike algebras, {\em Math. Z.} {\bf 240} (2002), 579--610.

\bibitem{FLO} O. Foda, B. Leclerc, M. Okado, J.-Y. Thibon and T.
Welsh, Branching functions of $A_{(n-1)}^{(1)}$ and Jantzen-Seitz
problem for Ariki-Koike algebras, {\em Adv. Math.} (2) {\bf 141} (1999),
322--365.

\bibitem{FSS} J. Fuchs, B. Schellekens and C. Schweigert, From Dynkin diagram symmetries to fixed
point structures, Comm. Math. Phys. {\bf 180}, 39--97 (1996).

\bibitem{Ge} M. Geck, On the representation theory of Iwahori-Hecke algebras
  of extended finite Weyl groups, {\em Represent. Theory} {\bf 4} (2000), 370--397.

\bibitem{Gen} G. Genet, On decomposition matrices for graded algebras, {\em J. Alg.} (1) {\bf 274}
(2004), 523--542.

\bibitem{GJ} G. Genet and N. Jacon, Modular representations of
cyclotomic Hecke algebras of type $G(r,p,n)$, Inter. Math. Res. Notices, (2006), 1--18.

\bibitem{GL} J. J. Graham and G. I. Lehrer, Cellular algebras, {\em Invent. Math.} {\bf 123} (1996), 1--34.

\bibitem{Hu1} J. Hu, A Morita equivalence theorem for Hecke algebras of type $D_n$ when
$n$ is even, {\em Manuscr. Math.} {\bf 108} (2002), 409--430.

\bibitem{Hu2} J. Hu, Crystal basis and simple modules for Hecke algebra of
  type $D_n$, {\em J. Alg.} (1) {\bf 267} (2003), 7--20.

\bibitem{Hu3} J. Hu, Modular representations of Hecke algebras of type
  $G(p,p,n)$, {\em J. Alg.} (2) {\bf 274} (2004), 446--490.

\bibitem{Hu4} J. Hu, Branching rules for Hecke Algebras of Type
$D_{n}$, {\em Math. Nachr.} {\bf 280} (2007), 93--104.

\bibitem{Hu5} J. Hu, Crystal basis and simple modules for Hecke algebra of
  type $G(p,p,n)$, {\em Representation Theory} {\bf 11} (2007), 16--44.

\bibitem{Hu6} J. Hu, The representation theory of the cyclotomic Hecke algebras of
  type $G(r,p,n)$, to appear in a Proceeding volume on the International Conference and Instructional workshop on
  discrete groups (Morningside center of Mathematics, Beijing, July 17-August 4, 2006).

\bibitem{Ja0} N. Jacon, {\em Repr\'esentations modulaires des alg\`ebres
de Hecke et des alg\`ebres de Ariki-Koike}, Ph.D. thesis,
Universit\'e Claude Bernard Lyon I, 2004.

\bibitem{Ja} N. Jacon, On the parameterization of the simple modules for
Ariki-Koike algebras at roots of unity, {\em J. Math. Kyoto Univ.}, {\bf
44} (2004), 729--767.

\bibitem{Ja1} N. Jacon, private email communication.

\bibitem{LLT} A. Lascoux, B. Leclerc and J.-Y. Thibon, Hecke
algebras at roots of unity and crystal bases quantum affine
algebras, Comm. Math. Phys. {\bf 181}, 205--263 (1996).

\bibitem{MM} T. Misra and K.C. Miwa, Crystal bases for the basic representations of $U_q(\ksl_n)$,
 Comm. Math. Phys. {\bf 134}, 79--88 (1990).

\bibitem{NS2} S. Naito and D. Sagaki, Lakshmibai-Seshadri paths fixed by a diagram automorphism,
  {\em J. Alg.} {\bf 245} (2001), 395--412.

\bibitem{NS1} S. Naito and D. Sagaki, Standard paths and standard monomials fixed by a diagram automorphism,
  {\em J. Alg.} {\bf 251} (2002), 461--474.

\bibitem{P} C. Pallikaros, Representations of Hecke algebras of type $D_n$,
  J. Alg. {\bf 169}, 20--48 (1994).

\bibitem{V} M. Vazirani, Irreducible modules over the affine Hecke algebra:
  a strong multiplicity one result, Ph.D. thesis, University of California at
  Berkeley, 1998.

\end{thebibliography}

\begin{thebibliography}{99}

\bibitem[{[}G{]}]{Gen1} G. Genet, On decomposition matrices for graded algebras, {\em J. Alg.} (1) {\bf 274}
(2004), 523--542.

\end{thebibliography}

\bigskip\bigskip

\vfill
\noindent
Department of Applied Mathematics,\\[2pt]
\noindent
Beijing Institute of Technology\\[2pt]
\noindent
Beijing, 100081\\[2pt]
The People's Republic of China\\[3pt]
\noindent
E-mail: xiaoyi.cui@gmail.com

\end{document}